\documentclass[]{amsart}

 \usepackage{amssymb}

 \usepackage{url}

\newtheorem{theorem}{Theorem}[section]
\newtheorem{corollary}[theorem]{Corollary}

\newtheorem{proposition}[theorem]{Proposition}

\theoremstyle{remark}

\theoremstyle{definition}

\numberwithin{equation}{section}
\numberwithin{theorem}{section}

\theoremstyle{plain}

\begin{document}

%%%%%%%%%%%%%%%%%%%%%%%%%
% Subject classification 
%%%%%%%%%%%%%%%%%%%%%%%%%

% Provide an AMS subject classification with one or two primary classification 
% numbers and, optionally, one or more secondary classification numbers. 
% Use the following format:  "Primary 42B25. Secondary 42B60, 20E26"


\date{Preprint September 1, 2011}
%%%%%%%%%
% Title
%%%%%%%%%

% Title, in lower case, with no explicit linebreaks (\\).  If the title
% is too long to be used as a running head, add a short version of the
% title in brackets, as in \title[shorttitle]{fulltitle}.

\title{Quadrature rules with (not too many) derivatives}

%%%%%%%%%%%%%%%%%%%%%%%%%%%%%%
% Author names and addresses 
%%%%%%%%%%%%%%%%%%%%%%%%%%%%%%

% Provide one separate \author{...} \address{...} \email{....} entry for each
% author, i.e., do not combine multiple authors.  Separate address lines by double
% slashes.  Do not attach footnotes to author  names. (For acknowledgements use
% the "\thanks" construct below.)
%
\author{Matthew Wiersma}
\address{Department of Pure Mathematics\\
University of Waterloo\\
Waterloo, Ontario, Canada N2L 3G1}
\email{mwiersma@uwaterloo.ca}
\thanks{The author wrote this paper under the supervision of Erik Talvila at the University of the Fraser Valley while supported by an Undergraduate Student Research Award from the Natural Sciences and Engineering Research Council of Canada.}

%%%%%%%%%%%%%%%%%%%%
% Acknowledgements
%%%%%%%%%%%%%%%%%%%

% Use \thanks for acknowledgements as footnotes to the title page.  
% (Note that footnotes inside \author or \title macros are not
% allowed.)
%
% In case of multiple author papers, phrase the acknowledgement to 
% say "The first author was supported by ...  The second author was
% supported by ..."

%%%%%%%%%%%%%
% Abstract 
%%%%%%%%%%%%%
%
% Abstracts should not contain macros (so that they can be processed independently
% of the paper.) Avoid displayed math and references in the abstract.

\begin{abstract}
Quadrature formulas for $\int_a^b f(x)\,dx$ where derivative terms need only be evaluated at $a$ and $b$ in the composite rule are identified. Error bounds are given when $f\!:[a,b]\to\mathbb{R}$ satisfies $f^{(n-1)}$ is absolutely continuous so that $f^{(n)}\in L^p([a,b])$, and when $f^{(n-1)}$ is merely continuous.
\end{abstract}

\maketitle

% New definition of square root:
% it renames \sqrt as \oldsqrt
\let\oldsqrt\sqrt
% it defines the new \sqrt in terms of the old one
\def\sqrt{\mathpalette\DHLhksqrt}
\def\DHLhksqrt#1#2{%
\setbox0=\hbox{$#1\oldsqrt{#2\,}$}\dimen0=\ht0
\advance\dimen0-0.2\ht0
\setbox2=\hbox{\vrule height\ht0 depth -\dimen0}%
{\box0\lower0.4pt\box2}}

\newcounter{projects}
\newcommand{\fn}{\!:\!}
\newcommand{\bv}{{\mathcal BV}}
\newcommand{\be}{\begin{equation}}
\newcommand{\ee}{\end{equation}}
\providecommand{\abs}[1]{\lvert#1\rvert}
\providecommand{\norm}[1]{\lVert#1\rVert}
\newcommand{\intab}{\int_a^b}
\newcommand{\polyk}{{\mathcal P}_k}
\newcommand{\poly}[1]{{\mathcal P}_{#1}}
\newcommand{\R}{{\mathbb R}}

\section{Introduction}

We are interested in finding numerical integration schemes for $\int_a^b f(x)\,dx$ where derivative terms in the composite rule need only be evaluated at the endpoints, $a$ and $b$. Let $p$ be a polynomial of degree $n$ with leading coefficient $1/n!$ and $f\!:[a,b]\to\mathbb{R}$ such that $f^{(n-1)}$ is absolutely continuous. Repeated integration by parts shows that
\begin{eqnarray}
\int_a^b f(x)\, dx &=& (b-a)p^{(n-1)}(1)f(b)-(b-a)p^{(n-1)}(0)f(a) \notag \\
&& -(b-a)^2 p^{(n-2)}(1)f'(b)+(b-a)^2 p^{(n-2)}(0)f'(a)\\
& & +\ldots -(-1)^n (b-a)^n p(1)f^{(n-1)}(b) \notag \\
& & +(-1)^n (b-a)^n p(0)f^{(n-1)}(a)+E_{(a,b)}(p,f) \notag
\end{eqnarray}
where the {\it error term} is
$$E_{(a,b)}(p,f)=(-1)^n(b-a)^{n+1}\int_0^1 p(x)f^{(n)}\big(a+(b-a)x\big)\, dx.$$
This method is used in \cite{cruzuribe}. The composite rule is given in equation (3.1). Taking $p(x)=x(x-1)/2$ in equation (1.1) gives the trapezoidal rule. The Euler-Maclaurin formula follows from taking $p(x)=B_n(x)/n!$ or $p(x)=(B_n(x)-B_n)/n!$ where $B_n(x)$ are Bernoulli polynomials and $B_n$ are Bernoulli numbers. The midpoint rule and Simpson's rule also follow from equation (1.1) \cite{trap}.

We find the integration schemes where the derivative terms need only be evaluated at the endpoints of the interval $[a,b]$ in the composite rule to be similar to the Euler-Maclaurin formula (Corollary 3.2). The most important of these integration schemes are essentially the Euler-Maclaurin formula. Similar formulas are derived in \cite{dedic}.

Section 2 outlines necessary facts on Bernoulli polynomials. The desired formulas are found in Section 3. In section 4 we give sharp error bounds when $f^{(n)}\in L^p([a,b])$ for $1\leq p\leq\infty$ and find asymptotic estimates for these bounds as $n\to\infty$. In section 5 we weaken the assumption that $f^{(n-1)}$ is absolutely continuous to merely assuming that $f^{(n-1)}$ is continuous. Then $f^{(n)}$ exists as a distribution and the error term can be evaluated with the distributional Denjoy integral. The error is estimated in terms of the Alexiewicz norm.

\section{Bernoulli Polynomials}

The Bernoulli polynomials can be defined recursively by $B_0(x)=1$ and
\begin{equation}
\frac{d}{dx}B_n(x)=nB_{n-1}(x)
\end{equation}
such that
\begin{equation}
\int_0^1 B_n(x)\,dx=0
\end{equation}
for $n\geq 1$. The Bernoulli numbers are defined as
$$ B_n\equiv B_n(0).$$
It then follows that
\begin{equation}
B_n(x)=\sum_{k=0}^n \binom{n}{k} B_k x^{n-k}.
\end{equation}
Integrating both sides of equation (2.3) on the interval $[0,1]$ leads to the equation
\begin{equation}
\sum_{k=0}^n \binom{n+1}{k} B_k=0
\end{equation}
when $n\geq 1$. The Bernoulli numbers may then be calculated recursively with equation (2.4).

It follows from induction that
\begin{equation}
B_n(1-x)=(-1)^nB_n(x).
\end{equation}
We may then conclude from equations (2.1), (2.2), and (2.5) that $B_n=0$ when $n>1$ is odd. It also follows from induction that $B_n(x)$ is monotonic on $[0,1/2]$ when $n$ is even, and $B_n(x)$ does not change sign on $[0,1/2]$ when $n$ is odd.

Other properties of Bernoulli polynomials include
\begin{eqnarray}
|B_{2n}(x)| &\leq & |B_{2n}| \textnormal{ for } x\in[0,1]\\
|B_{2n}(x)-B_{2n}| &\leq& (2-2^{1-2n})|B_{2n}| \textnormal{ for } x\in[0,1]\\
\int_0^{1/2}B_n(x)\,dx &=& \frac{1-2^{n+1}}{2^n}\frac{B_{n+1}}{n+1}\\
\int_0^1 B_n(x)B_m(x)\,dx &=& \frac{(-1)^{n-1}m!n!}{(m+n)!}B_{m+n} \textnormal{ for } m,n\geq 1\\
(-1)^{n-1}\frac{(2\pi)^{2n}B_{2n}(x)}{2(2n)!} &\to& \cos(2\pi x) \textnormal{ as } n\to\infty\\
(-1)^{n-1}\frac{(2\pi)^{2n+1}B_{2n+1}(x)}{2(2n+1)!} &\to& \sin(2\pi x) \textnormal{ as } n\to\infty.
\end{eqnarray}
The convergence in equations (2.10) and (2.11) is uniform on any compact set. All these results may be found in \cite{nist}.

\section{Main Theorem}

It is useful to define some new terms before proceeding further. Let $p$ be a polynomial of degree $n$ with leading coefficient $1/n!$ and $f\!:[a,b]\to\mathbb{R}$ such that  $f^{(n-1)}$ is absolutely continuous. Then
\begin{eqnarray*}
I_{(a,b)}(p,f) &:=& (b-a)p^{(n-1)}(1)f(b)-(b-a)p^{(n-1)}(0)f(a) \notag \\
&& -(b-a)^2 p^{(n-2)}(1)f'(b)+(b-a)^2 p^{(n-2)}(0)f'(a) \notag \\
& & +\ldots -(-1)^n (b-a)^n p(1)f^{(n-1)}(b)\notag \\
&& +(-1)^n (b-a)^n p(0)f^{(n-1)}(a)\notag
\end{eqnarray*}
is said to be the {\it integration scheme generated by $p$}. We may now phrase equation (1.1) as
$$\int_a^b f(x)\,dx = I_{(a,b)}(p,f)+E_{(a,b)}(p,f).$$
Let $h=(b-a)/N$ for natural number $N$. The {\it composite rule} is
\begin{equation}
\int_a^b f(x)\,dx=I_{(a,b)}^N(p,f)+E_{(a,b)}^N(p,f)
\end{equation}
where
$$I_{(a,b)}^N(p,f)=\sum_{k=0}^{N-1}I_{(a+kh,a+(k+1)h)}(p,f)$$
and
$$E_{(a,b)}^N(p,f)=\sum_{k=0}^{N-1}E_{(a+kh,a+(k+1)h)}(p,f).$$

In general, the derivative terms in the composite rule need only be evaluated at the endpoints, $a$ and $b$, if and only if $p^{(\ell)}(0)=p^{(\ell)}(1)$ for $0\leq \ell\leq n-2$. In other words, this happens when the derivative terms of $f$ form a telescoping series in the composite rule. This motivates us to say that $p$ is a {\it telescoping polynomial} for $I_{(a,b)}(p,f)$ when $p^{(\ell)}(0)=p^{(\ell)}(1)$ for $0\leq \ell\leq n-2$. In general, the function values of $f$ in the composite rule cannot form a telescoping series since $p^{(n-1)}$ is linear. In the following theorem we charaterise the telescoping polynomials.

\begin{theorem}
The polynomials
$$\frac{B_n(x)}{n!}+c,$$
where $c$ is an arbitrary constant, form the set of degree $n$ telescoping polynomials.
\end{theorem}

\begin{proof}
Let $p$ be a polynomial of degree $n$ with leading coefficient $1/n!$. Write
$$p(x)=\sum_{k=0}^n \frac{\alpha_{k}}{(n-k)!k!}x^{n-k}$$
where $\alpha_0=1$. Then
$$p^{(\ell)}(x)=\sum_{k=0}^{n-\ell}\frac{\alpha_{k}}{(n-\ell-k)!k!}x^{n-\ell-k}.$$
Equating $p^\ell(0)$ and $p^\ell(1)$, it follows that $p$ is a telescoping polynomial if and only if 
$$\sum_{k=0}^{n-\ell-1}\frac{\alpha_{k}}{(n-\ell-k)!k!}=0$$
for $0\leq \ell\leq n-2$, i.e., $1\leq n-\ell-1\leq n-1$. The $\alpha_k$ are then defined recursively as in equation (2.4) for $0\leq k\leq n-1$ and $\alpha_n$ is arbitrary.
\end{proof}

The composite rules produced by these polynomials are found by evaluating equation (3.1) with $p(x)=B_n(x)/n!+c$.

\begin{corollary}
Let $f:[a,b]\to\mathbb{R}$ such that $f^{(n-1)}$ is absolutely continuous. Then for each constant $c$,
\begin{eqnarray}
\int_a^b f(x)\,dx &=& \left(\frac{b-a}{2N}\right)\big(f(a)+f(b)\big)+\left(\frac{b-a}{N}\right)\sum_{k=1}^{N-1}  f(a+kh) \nonumber \\
&&+\sum_{k=2}^{n}\frac{B_k}{k!}\left(\frac{b-a}{N}\right)^{k}\big(f^{(k-1)}(a)-f^{(k-1)}(b)\big)  \\
&&+c(b-a)^n\big(f^{(n-1)}(a)+f^{(n-1)}(b)\big)+E_{(a,b-a)}^N(B_n(x)/n!+c,f). \nonumber
\end{eqnarray}
\end{corollary}
The fact that $c$ is arbitrary begs the question of how $c$ should be chosen. We will focus on the cases when $c=-B_n/n!$ ($n\geq 2$) and when $c=0$. The $f^{(n-1)}$ terms in equation (3.2) are not evaluated when $c=-B_n/n!$. In this way, taking $c=-B_n/n!$ gives the integration scheme with the least number of derivative terms generated by a degree $n$ polynomial. The integration scheme generated by $B_n(x)/n!$ is also generated by a polynomial of higher degree. This is not true when $c\neq 0$. Define $p_n(x)=B_n(x)/n!$ and $q_n(x)=B_n(x)/n!-B_n/n!$. These polynomials relate to the choices of $c$ indicated above. When $n\geq 3$ is odd, $p_n=q_n$. Moreover, the integration schemes generated by $p_{2n}$, $p_{2n+1}$ and $q_{2n+2}$ are equal, and are essentially the Euler-Maclaurin formula. The following proposition shows some cases of when $p_n$ minimises the error estimate in equation (4.1).

\begin{proposition}
(a) Choosing $c=0$ minimises $\|p_n+c\|_2$ for each $n\geq 0$.\\
(b) Moreover, $c=0$ minimises $\|p_n+c\|_r$ for each odd $n$, $1\leq r\leq \infty$.
\end{proposition}
\begin{proof}
(a) It follows from equation (2.2) that
$$\|p_n+c\|_2^2=\int_0^1 [p_n(x)]^2\,dx+c^2.$$
The desired conclusion is then easily observed.

(b) The second assertion clearly follows from equation (2.5) when $r=\infty$. When $r<\infty$,
\begin{eqnarray*}
\|p_n+c\|_r^r &=& \int_{0}^{1/2}\left|p_n(x)+c\right|^r\,dx+\int_{1/2}^1 \left|p_n(x)+c\right|^r\,dx\\
&=& \int_{0}^{1/2}\left|p_n(x)+c\right|^r\,dx+\int_{0}^{1/2}\left|p_n(x)-c\right|^r\,dx
\end{eqnarray*}
due to the asymmetry of $p_n$ about $1/2$ when $n$ is odd. The minimum of \linebreak $\abs{p_n(x)+c}^r+\abs{p_n(x)-c}^r$ occurs at $c=0$. Hence,
$$ \norm{p_n}_r^r = 2\int_0^{1/2} |p_n(x)|^r\,dx \leq \norm{p_n+c}_r^r.$$
The inequality becomes strict for nonzero $c$ when $r>1$.
\end{proof}

\section{$L^p$ Error Estimates}

We may use H\"{o}lder's inequality to estimate the error terms under the assumed conditions. Let $1/r+1/s=1$ for positive $r,s$. Choose $f$ so that $f^{(n-1)}$ is absolutely continuous and $f^{(n)}\in L^s([a,b])$. Then
\begin{equation}
\left|E_{(a,b)}(p,f)\right|\leq (b-a)^{n+1/r}\|p\|_r\|f^{(n)}\|_s
\end{equation}
where $\|p\|_r=\left(\int_0^1|p(x)|^r\,dx\right)^{1/r}$ if $r<\infty$, $\|f^{(n)}\|_s=\left(\int_a^b|f^{(n)}(x)|^s\,dx\right)^{1/s}$ if $s<\infty$, $\|p\|_\infty=\max_{x\in [0,1]}|p(x)|$, $\|f^{(n)}\|_\infty=\sup_{x\in[a,b]}|f^{(n)}(x)|$. Proposition 4.1 shows that this bound is sharp.

H\"{o}lder's inequality for series shows that
$$\sum_{k=0}^{N-1}\|f^{(n)}\chi_{[a+kh,a+(k+1)h]}\|_s\leq N^{1/r}\|f^{(n)}\|_s$$
when $h=(b-a)/N$. We may then estimate the error in the composite rules as
\begin{equation}
\left|E_{(a,b)}^N(p,f)\right|\leq \frac{(b-a)^{n+1/r}}{N^{n}}\|p\|_r\|f^{(n)}\|_s.
\end{equation}

\begin{proposition}
The bound given in equation (4.1) is sharp.
\end{proposition}

\begin{proof}
When $r=1$, equality holds when $f$ satisfies $f^{(n)}(a+(b-a)x)=\textnormal{sgn}(p(x))$ for $x\in [0,1]$.
When $1<r<\infty$, choose $f$ so that $f^{(n)}(a+(b-a)x)=|p(x)|^{r/s}\textnormal{sgn}(p(x))$. Then equality holds in equation (4.1) \cite[Theorem 2.3]{liebloss}. Let $\{\psi_k\}$ be the delta sequence given by $\psi_k(x)=1/k$ for $|x|<1/(2k)$ and $\psi_k(x)=0$ otherwise, and let $x_0\in [0,1]$ maximise $|p|$ over $[0,1]$. Choose $f_k$ to satisfy
$f_k^{(n)}(a+(b-a)x)=\psi_k(x-x_0)$ if $x_0\in (0,1)$ and $f_k^{(n)}(a+(b-a)x)=2\psi_k(x-x_0)$ otherwise. Then $\|f_k^{(n)}\|_1=(b-a)$ and $\int_0^1 p(x)f_k^{(n)}(x)\,dx\to p(x_0)$ as $k\to\infty$.
\end{proof}

In general, calculating $\|p\|_r$ is a difficult problem. The following proposition gives some cases of when we are able to calculate this exactly for $p_n$ and $q_n$.

\begin{proposition}
(a) $\|q_n\|_1=|B_{n}|/n!$ when $n$ is even.\\
(b) $\|p_n\|_1=(1-2^{n+1})|B_{n+1}|/((n+1)2^{n-1})$ when $n$ is odd.\\
(c) $\|p_n\|_2=\sqrt{|B_{2n}|/(2n)!}$ and $\|q_n\|_2=\sqrt{|B_{2n}|/(2n)!+(B_{n}/n!)^2}$.\\
(d) When $n$ is even, $\|p_n\|_\infty=|B_n|/n!$ and $\|q_n\|_\infty=|q_n(1/2)|\leq (2-2^{1-n})|B_n|/n!$.
\end{proposition}
\begin{proof}
(a) Let $n$ be even. Then $q_n$ does not change sign on $[0,1]$ by equation (2.6). Hence,
$$\int_0^1 |q_n(x)|\,dx = \left| \int_0^1 q_n(x)\,dx \right| = |B_n|/n!.$$

(b) Let $n$ be odd. Then $p_n$ is asymmetric about 1/2 while not changing sign on $[0,1/2]$. Hence,
$$\int_0^1 |p_n(x)|\,dx = 2\left| \int_0^{1/2} q_n(x)\,dx \right| = \frac{1-2^{n+1}}{2^{n-1}} \frac{|B_{n+1}|}{n+1}.$$
The final equality follows from equation (2.8).

(c) The first inequality is evident from equation (2.9). The second then follows from the first while noting equation (2.2).

(d) The first equality is essentially equation (2.6). The second equality comes from $q_n$ being monotonic over $[0,1/2]$ while being symmetric about 1/2 with $q_n(0)=0$. The inequality is essentially equation (2.7).
\end{proof}

We now shift our attention to finding asymptotic estimates of $\|p_n\|_r$ and $\|q_n\|_r$ as $n\to\infty$. Since the convergence in equations (2.10) and (2.11) is uniform,
\begin{eqnarray}
\int_0^1 |p_{2n}(x)|^r\,dx &\sim & \left(\frac{2}{(2\pi)^{2n}}\right)^r \int_0^1 |\cos(2\pi x)|^r\,dx \textnormal{ as } n\to\infty \\
\int_0^1 |p_{2n+1}(x)|^r\,dx &\sim & \left(\frac{2}{(2\pi)^{2n+1}}\right)^r \int_0^1 |\sin(2\pi x)|^r\,dx \textnormal{ as } n\to\infty\\
\int_0^1 |q_{2n}(x)|^r\,dx &\sim & \left(\frac{2}{(2\pi)^{2n}}\right)^r \int_0^1 |\cos(2\pi x)+1|^r\,dx \textnormal{ as } n\to\infty
\end{eqnarray}
and
\begin{eqnarray}
\|p_n\|_\infty &\sim& \frac{2}{(2\pi)^n} \textnormal{ as } n\to\infty \\
\|q_{2n}\|_\infty &\sim& \frac{4}{(2\pi)^{2n}} \textnormal{ as } n\to\infty.
\end{eqnarray}
Evaluating equations (4.3), (4.4) and (4.5) at $r=1$ gives:
\begin{eqnarray}
\|p_n\|_1 &\sim& \frac{8}{(2\pi)^{n+1}} \textnormal{ as } n\to\infty \\
\|q_{2n}\|_1 &\sim& \frac{2}{(2\pi)^{2n}} \textnormal{ as } n\to\infty .
\end{eqnarray}
It is well known that
\begin{equation}
 \int_0^{\pi/2}\sin^r x\,dx=\int_0^{\pi/2}\cos^r x\,dx=\frac{1\cdot 3\cdot 5\cdots (r-1)}{2\cdot 4\cdot 6\cdots r}\frac{\pi}{2}
\end{equation}
when $r\geq 2$ is even and
\begin{equation}
\int_0^{\pi/2}\sin^r x\,dx=\int_0^{\pi/2}\cos^r x\,dx=\frac{2\cdot 4\cdot 6\cdots (r-1)}{1\cdot 3\cdot 5\cdots r}
\end{equation}
when $r\geq 3$ is odd. Using equations (4.10) and (4.11) we find the following asymptotic estimates for $\|p_n\|_r$ and $\|q_n\|_r$ when $r\geq 2$ is an integer:
\begin{equation}
\|p_n\|_{r} \sim \frac{2}{(2\pi)^{n}} \left(\frac{1\cdot 3\cdot 5\cdots (r-1)}{2\cdot 4\cdot 6\cdots r}\right)^{1/r} \textnormal{ as } n\to\infty
\end{equation}
when $r\geq 2$ is even, and
\begin{equation}
\|p_n\|_{2r+1} \sim \frac{2}{(2\pi)^{n}} \left(\frac{2\cdot 4\cdot 6\cdots (r-1)}{1\cdot 3\cdot 5\cdots r}\frac{2}{\pi}\right)^{1/r} \textnormal{ as } n\to\infty 
\end{equation}
when $r\geq 3$ is odd.

\section{Generalised Error Bounds}

We initially assumed that $f^{(n-1)}$ was absolutely continuous on $[a,b]$. We may weaken this assumption to merely assuming $f^{(n-1)}$ is continuous if we evaluate integrals in the distributional Denjoy sense. The distributional Denjoy integral allows integrating derivatives of any continuous functions where derivatives are taken in the distributional sense.

Define the set of {\it test functions} to be $\mathcal{D}([a,b])=\{\phi\in C^{\infty}([a,b])\mid \phi\textnormal{ has support in } (a,b) \}$ where the {\it support} of a function is the closure of the set of points where the function is nonzero. A sequence $\{\phi_n\}\subset \mathcal{D}([a,b])$ is said to converge to $\phi\in\mathcal{D}([a,b])$ when there exists a compact set $K\subset (a,b)$ so that $\phi_n$ has support in $K$ for each $n$ and $\phi_n^{(m)}\to \phi^{(m)}$ uniformly. The {\it distributions} are the continuous linear functionals on the space of test functions; the set of distributions is denoted $\mathcal{D}'([a,b])$. Let $f\!:[a,b]\to\mathbb{R}$ be an integrable function. Then $\langle f,\phi\rangle=\int_a^b f(x)\phi(x)\,dx$ defines a distribution.

Let $T$ be a distribution. Then $T'$ so that $\langle T',\phi\rangle=-\langle T,\phi'\rangle$ for each $\phi\in\mathcal{D}([a,b])$ is said to be the {\it distributional derivative} of $T$. If $T$ is taken to be an integrable function $f\!: [a,b]\to \mathbb{R}$ with a pointwise derivative at $x_0$ then that pointwise value may be identified by evaluating distribution $f$ with a delta sequence of test functions about $x_0$.

The definition of the {\it distributional Denjoy integral} is simple. If $f$ is the distributional derivative of a continuous function $F\!:[a,b]\to\mathbb{R}$ then $\int_a^b f=F(b)-F(a)$ defines the integral over a compact interval. Take $f$ to be the distributional derivative of a continuous function $F\!:[a,b]\to\mathbb{R}$ and let $g\!:[a,b]\to\mathbb{R}$ be of bounded variation. The integration by parts formula for the distributional Denjoy integral is given by
\begin{equation}
\int_a^b fg = F(b)g(b)-F(a)g(a)-\int_a^b F\,dg.
\end{equation}
In this way, equations (1.1) and (3.1) still hold with this weakened assumption that $f^{(n-1)}$ be merely continuous. The Beesack-Darst-Pollard inequality for the distributional Denjoy integral is given by
\begin{equation}
\left|\int_a^b fg\right|\leq \left|\int_a^b f\right|\inf_{[a,b]}|g|+\|f\|_{\mathcal{A}_c}V_{[a,b]}g.
\end{equation}
Here $\|\cdot\|_{\mathcal{A}_c}$ is the Alexiewicz norm given by
$$ \|f\|_{\mathcal{A}_c}=\max_{x,y\in [a,b]}\left|F(x)-F(y)\right| $$
and $V_{[a,b]}g$ is the total variation of $g$. All this on the distributional Denjoy integral is found in \cite{denjoy}.

Applying equation (5.2) to $E_{(a,b)}(p_n+c, f)$ gives the estimate
\begin{equation}
\left|E_{(a,b)}(p_n+c,f)\right| \leq (b-a)^n\|f^{(n)}\|_{\mathcal{A}_c}V_{[0,1]}p_n
\end{equation}
when $p_n+c$ has a real root in $[0,1]$. This is true for both $p_n$ and $q_n$. It then follows that
\begin{equation}
\left|E_{(a,b-a)}^N(p_n+c,f)\right|\leq \frac{(b-a)^n}{N^{n-1}}\|f^{(n)}\|_{\mathcal{A}_c}V_{[0,1]}p_n
\end{equation}
when $p_n+c$ has a real root in $[0,1]$.

The total variation of $p_n$ may be calculated from the maxima and minima of Bernoulli polynomials on $[0,1]$. See \cite{lehmer} for estimates of these upper and lower bounds.

We use the uniform convergence in equations (2.10) and (2.11) to give the asymptotic estimate
$$ V_{[0,1]}p_n \sim \frac{8}{(2\pi)^n} \textnormal{ as } n\to\infty. $$

\section*{Acknowledgements}
The author extends his deepest gratitude to Erik Talvila, without whom this paper would not be possible.


\begin{thebibliography}{99}

\bibitem{cruzuribe}
D. Cruz-Uribe and C.J. Neugebauer, {\it Sharp error bounds for the trapezoidal
rule and Simpson's rule}, JIPAM. J. Inequal. Pure Appl. Math. {\bf 3}(2002),
Article 49, 22~pp.

\bibitem{dedic}
Lj. Dedi\`{c}, M. Mati\`{c} and J. Pe\v{c}ari\`{c},
{\it On Euler trapezoid formulae},
Appl. Math. Comput. {\bf 123}(2001), 37--62.

\bibitem{nist}
K. Dilcher {\it Bernoulli and Euler polynomials}, in: NIST handbook of mathematical functions, U.S. Dept. Commerce, Washington, DC, 2010, pp. 587--599.

\bibitem{lehmer}
D.H. Lehmer,
{\it On the maxima and minima of Bernoulli polynomials},
Amer. Math. Monthly {\bf 47}(1940). 533--538.

\bibitem{liebloss}
E.H. Lieb and M. Loss,
{\it Analysis},
Providence, American Mathematical Society, 2001.

\bibitem{denjoy}
E. Talvila, {\it The distributional Denjoy integral}, Real Anal. Exchange {\bf 33}(2008), 51--82.

\bibitem{quadoptimal}
E. Talvila and M. Wiersma,
{\it Optimal error estimates for corrected trapezoidal rules},
(to appear).

\bibitem{trap}
E. Talvila and M. Wiersma,
{\it Simple proofs of basic quadrature formulas},
(to appear).
\end{thebibliography}
\end{document}